# Robust Integration of High-level Dispatchable Renewables in Power System Operation

Hongxing Ye, *Student Member, IEEE,* Jianhui Wang, *Senior Member, IEEE,* and Zuyi Li, *Senior Member, IEEE*

*Abstract*—The increasing penetration of Renewable Energy Sources (RES) requires more Flexibility Resources (FR), generally thermal units and storages, must be kept in the system to accommodate the uncertainties from RES. The challenge is how the system can survive when the RES level is very high. In this paper, RESs are considered as full-role market participants. They can bid in the day-ahead market, and the powers they deliver to the market are controllable up to their maximum available powers. Therefore, RESs are effectively dispatchable and can function as FR providers. To integrate dispatchable renewables, a two-stage robust Unit Commitment (UC) and dispatch model is established. In the first stage, a base UC and dispatch is determined. In the second stage, all FRs including RESs are used to accommodate the uncertainties, which is a Mixed-Integer Programming (MIP) problem. It is proved that the solution to the max-min problem can be identified directly whether the strong duality holds or not for the inner minimization problem. The solution robustness is guaranteed by including only one extra scenario. Numerical results show the effectiveness of the proposed model and its advantages over the traditional robust UC model with high level RES penetration.

## Nomenclature

**Indices**

| | |
|---|---|
| $t$ | index of time intervals |
| $i, l, m$ | index of thermal units, transmission lines, and buses |
| $r$ | index of RES units |
| $f$ | index of fast startup units |

**Functions and sets**

| | |
|---|---|
| $\mathcal{F}_1$ | feasible region for base-case UC and dispatch |
| $\mathcal{F}_2$ | feasible region for robust UC and dispatch |
| $\mathcal{F}_3$ | feasible region for feasibility check problem |
| $\mathcal{F}_4$ | feasible region for worst-case UC and dispatch |
| $\mathcal{U}$ | uncertainty set |
| $C_f^F(\cdot)$ | cost related with fast startup unit $f$ |
| $C_i^I(\cdot)$ | cost related to UC for thermal unit $i$ |
| $C_i^P(\cdot)$ | cost related to dispatch for thermal unit $i$ |
| $C_r^R(\cdot)$ | bid-based cost related to RES unit $r$ |
| $\mathcal{G}_m$ | set of thermal units located at bus $m$ |
| $\mathcal{Q}_m$ | set of fast startup units located at bus $m$ |
| $\mathcal{R}_m$ | set of RES units located at bus $m$ |

**Constants**

| | |
|---|---|
| $\alpha$ | confidence interval parameter |
| $\beta$ | RES energy level |
| $N_R$ | number of RES units |
| $N_T$ | number of time intervals |
| $d_{m,t}$ | aggregated equivalent load at bus $m$ time $t$ |
| $F_l$ | transmission line flow limit for line $l$ |
| $\Gamma_{l,m}$ | shift factor for line $l$ with respect to bus $m$ |
| $P_i^{\min}$ | minimum power output for unit $i$ |
| $P_i^{\max}$ | maximum power output for unit $i$ |
| $r_i^u, r_i^d$ | ramping-up/down limits between sequential intervals for unit $i$ |
| $R_i^u, R_i^d$ | ramping-up/down limits for uncertainty accommodation for unit $i$ |
| $u_{r,t}$ | uncertainty bound of the power output for unit $r$ at time $t$ |
| $\tilde{P}_{r,t}^R$ | expected power output for unit $r$ at time $t$ |

**Variables**

| | |
|---|---|
| $\epsilon_{r,t}$ | power output uncertainty for unit $r$ at time $t$ |
| $\boldsymbol{\epsilon}$ | a compactor vector form of all uncertainties |
| $\bar{P}_{r,t}$ | the uncertain maximum available power output for unit $r$ at time $t$ |
| $I_{i,t}$ | on/off status indicator for unit $i$ at time $t$ |
| $x_{i,t}^{\text{on}}, x_{i,t}^{\text{off}}$ | number of hours unit $i$ has been on/off at time $t$ |
| $y_{i,t}, z_{i,t}$ | start-up/shut-down indicators for unit $i$ at time $t$ |
| $P_{i,t}$ | generation dispatch for unit $i$ at time $t$ |
| $\hat{P}_{i,t}$ | re-dispatch for unit $i$ at time $t$ when uncertainty is revealed |
| $I_{f,t}^F$ | on/off status indicator for unit $f$ at time $t$ when uncertainty is revealed |
| $P_{f,t}^F$ | dispatch of unit $f$ at time $t$ when uncertainty is revealed |
| $P_{r,t}^R$ | delivered power for unit $r$ at time $t$ |
| $\hat{P}_{r,t}^R$ | re-dispatched power of unit $r$ at time $t$ |
| $P_{m,t}^{\text{inj}}$ | net power injection at bus $m$ time $t$ |
| $\hat{P}_{m,t}^{\text{inj}}$ | net power injection at bus $m$ time $t$ after re-dispatch |
| $s_{m,t}$ | slack variable at bus $m$ time $t$ |

## I. Introduction

The Renewable Energy Sources (RES) are dramatically increasing in recent years. They produce competitive energy that is low in production cost and free of carbon emission. However, RES power output is intermittent and volatile. They also introduce more uncertainties into the system operation and electric markets. In the U.S. Day-ahead Market (DAM), these uncertainties are mainly from the day ahead forecasting errors. The Security-Constraint Unit Commitment (SCUC) and Economic Dispatch (ED) problems in DAM rely on

H. Ye and Z. Li are with Illinois Institute of Technology, Chicago, IL 60616, USA. J. Wang is with Argonne National Lab, Argonne, IL 60439, USA. (e-mail: hye9@hawk.iit.edu; jianhui.wang@anl.gov; lizu@iit.edu).

the accurate system data forecasting for the next day. The uncertainties from RESs pose new challenges for the SCUC and ED problems in DAM, which have attracted extensive attentions in recent years [1], [2], [3], [4], [5], [6].

The SCUC and ED problem is a daily task for the Independent System Operators (ISO) and Regional Transmission Organization (RTO) [7], [8], [9], [10]. In the SCUC problem, the on/off schedule of each generating unit, also called Unit Commitment (UC), in the next day is determined. The objective is to find the least cost UC combinations and power output levels for all the units to supply the load while respecting the system-wide constraints as well as unit-wise constraints. The system-wide constraints may include the load demand balance, transmission capacity limit, and reserve requirement. The unit-wise constraints are normally composed of generation capacity limits, ramping rate limits, and minimum on/off time limits [10], [9], [11].

To accommodate the uncertainties from RES, both scenario-based stochastic optimization and robust optimization are studied extensively for the SCUC problems [12], [3], [5], [4], [13], [14], [15], [16]. The main idea of the scenario-based stochastic SCUC is to generate a large number of sample points based on the Probability Density Function (PDF) of the random variables (i.e. uncertainties), then reduce the sample set to a much smaller one, which can be modeled and solved by modern Mixed-Integer Programming (MIP) solvers. In general, there are two main drawbacks of these scenario-based stochastic SCUC. The first one is that the PDF information is required, which is often unavailable in practice. The second one is that the scenario reduction may cause survivability issue of the system in some cases. Therefore, researchers in [5], [4] introduce the two-stage Robust SCUC (RSCUC) to overcome these drawbacks. In the traditional two-stage RSCUC, the UC solution is determined in the first stage, which leads to the least cost in the worst case. The ED is considered in the second stage, and it can be adjusted to accommodate any uncertainty in a pre-determined set. Because of the robustness associated with the solution, the two-stage RSCUC becomes another research focus in recent years. However, as the worst case is optimized, it has the issue of the over-conservativeness. In order to get an economical effective solution, authors in [14] introduce a unified approach to combine the scenario-based stochastic SCUC and RSCUC. [17] also introduce a new concept, the recourse cost requirement, to optimize the base case. Both the scenario-based stochastic SCUC and robust SCUC are computation intensive. Researchers in [13] reported the Column Generation (CG) algorithm to accelerate the RSCUC solution approach. On the other side, Affine Policy (AP) [18] is introduced to simplify the recourse actions in the robust approaches [15], [16], [19]. With strong assumptions, the intractable robust problem is converted to a convex and tractable one in the AP approach. However, the price of using AP is the reduction of the recourse actions, which often leads to non-zero optimality gap.

In most of the existing approaches, all the RES outputs (or uncertainties) must be accommodated. When the penetration level is low, RES generally can lower the total operation cost. However, it may not be true when the RES penetration level is high and the system has to accommodate all available energy from RES. A fundamental reason is that deliverable ramping capabilities must be kept in the system to accommodate the uncertainties. These ramping capabilities, also called flexibilities or reserves, are expensive when the uncertainty level is high.

A potential solution to this problem is to make RES dispatchable. That is, the RES power delivered to the system is controllable from zero to the maximum available power. As presented in [20], the wind power output could be controlled with the latest technologies. It is also possible to control the output of photovoltaic (PV) array by adjusting the array angle. On the other hand, too much zero-cost generation from RES could depress the Locational Marginal Price (LMP) at the injection point to a very low level or even zero or negative, which would ultimately affect the profit of RES. Thus, rather than acting as price takers, RESs could bid into the electricity markets as dispatchable resources. Recently, some researchers have explored the dispatchability of wind power [21], [22]. When RES is controllable in RSCUC, authors in [22] show that the strong duality can help to reformulate the second-stage problem.

In this paper, RESs are modeled as dispatchable resources with bid offers. A two-stage adaptive RSCUC is formulated. In the first stage, a set of RSCUC and ED solution is obtained, which optimally determines the least cost flexibility (or reserves) in the system. In the second stage, these flexibilities can accommodate any uncertainties within the confidence interval. The flexibility includes fast startup units, thermal units with available ramping capabilities, and dispatchable RESs. The main contributions of this paper are

1) A novel robust integration model of high-level dispatchable RESs is proposed. Fast startup units are considered in the second stage of RSCUC, which involves integer variables. It is proved that the solution to the non-convex max-min problem in the second stage can be obtained directly. Different from [22], the conclusion holds even without the strong duality of the inner-level minimization problem.
2) A fast solution approach is presented to solve the proposed RSCUC. Only two scenarios need to be considered. One is the base-case scenario and the other is the worst-case scenario. The robustness is guaranteed by modeling the worst-case scenario. The computational challenge in existing RSCUC is addressed by solving a converted single-level MIP problem with dispatchable RES.
3) It is revealed that the dispatchable RESs are more economical in terms of total system cost in the RSCUC. The bid offers of the RESs are modeled. The cost for accommodating all power from RESs is high when the uncertainty is large. With high level penetration, dispatching the RESs is potentially a new direction in the electricity markets.

The rest of this paper is organized as follows. In Section II, the new robust integration model for dispatchable RES is presented. Then, a solution approach is presented in Section III.

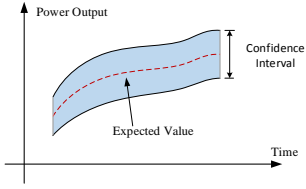

Fig. 1. The confidence interval and the expected value of RES power output

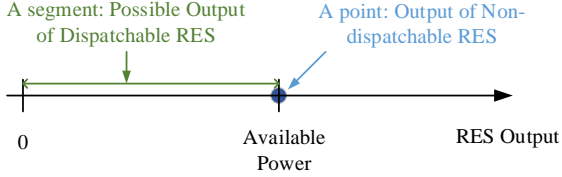

Fig. 2. The actual RES output delivered to the grid

Case studies are presented in Section IV. Finally, Section V concludes this paper.

## II. ROBUST INTEGRATION OF DISPATCHABLE RES

In the propose model, the RESs are treated as market participants with full roles. They are different from the traditional RSCUC [4], [5], [23], [17], where power outputs of RESs are considered as negative loads, and the bids of these power outputs are zeros. In this paper, the RESs can bid, and the output $P_{r,t}^R$ of RES unit $r$ at time $t$ is a decision variable. Due to the forecast errors, its upper bound $\bar{P}_{r,t}^R$ is an uncertain parameter.

$$\bar{P}_{r,t}^R = \tilde{P}_{r,t}^R + \epsilon_{r,t}, \forall r,t, \quad (1)$$

where $\tilde{P}_{r,t}^R$ is the expected power output for RES unit $r$ at $t$, and $\epsilon_{r,t}$ is the uncertainty. The uncertainty set $\mathcal{U}$ is defined as

$$\mathcal{U} := \{\boldsymbol{\epsilon} \in \mathbb{R}^{N_R N_T} : -u_{r,t} \leq \epsilon_{r,t} \leq u_{r,t}, \forall r,t\}.$$

Fig. 1 illustrates the possible RES power output. The uncertainty always falls within the confidence interval defined by $u_{r,t}$. The expected RES output is not modeled as the inelastic negative load anymore [17]. Instead, it is set as the upper bound of $P_{r,t}^R$ in the first stage

$$0 \leq P_{r,t}^R \leq \tilde{P}_{r,t}^R, \forall r,t. \quad (2)$$

In this paper, only intervals of the uncertainty are considered, which is similar to the interval optimization [24]. However, the largest power flow for each line is considered simultaneously. Hence, it is more optimistic than the interval optimization. The budget constraint in some robust approaches [5] is ignored. At the first glance, the uncertainty set in this paper is larger than that in other models, which may lead to the issue of over-conservativeness. However, the dispatchability of RESs reduces the conservativeness, especially when the RES penetration level is high. As shown in Fig. 2, power output of RES is inelastic in traditional robust approaches, and it is a point. When the RES is dispatchable, the RES power delivered to the grid is within a segment. Therefore, the feasible region of the re-dispatch problem is increased.

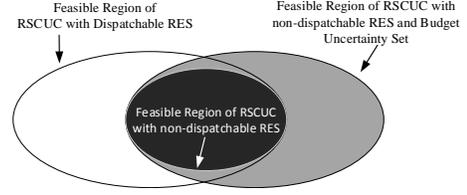

Fig. 3. Feasible region comparison for RSCUC

Fig. 3 compares the feasible regions for different RSCUC models. In robust approaches with non-dispatchable RES, adding the budget parameter (i.e. reducing the uncertainty set) can enlarge the feasible region. In comparison, dispatching RES also enlarges the feasible region.

The new RSCUC can be formulated as

$$\text{(P1):} \quad \min \quad \left\{ \begin{array}{l} \sum_i \sum_t \left( C_i^P(P_{i,t}) + C_i^I(I_{i,t}) \right) \\ + \sum_r \sum_t C_r^R(P_{i,t}^R) \end{array} \right\} \quad (3)$$

$$\text{s.t.} \quad \left\{ P_{i,t}, I_{i,t}, P_{i,t}^R \right\} \in \mathcal{F}_1 \cap \mathcal{F}_2 \quad (4)$$

where $\mathcal{F}_1$ is the feasible region for the base case in the first stage, and $\mathcal{F}_2$ guarantees the existence of the re-dispatch when the uncertainty is revealed in the second stage. The objective function (3) is to minimize the base cost. The solution to (P1) must be robust against any realization of the uncertainty in the second stage. $\mathcal{F}_1$ is defined as

$$\mathcal{F}_1 := \Big\{ (P, I, P^R) :$$

$$\sum_i P_{i,t} + \sum_r P_{r,t}^R = \sum_m d_{m,t}, \forall t \quad (5)$$

$$P_{m,t}^{\text{inj}} = \sum_{i \in \mathcal{G}_m} P_{i,t} + \sum_{r \in \mathcal{R}_m} P_{r,t}^R - d_{m,t}, \forall m,t \quad (6)$$

$$-F_l \leq \sum_m \Gamma_{l,m} P_{m,t}^{\text{inj}} \leq F_l, \forall l,t \quad (7)$$

$$I_{i,t} P_i^{\min} \leq P_{i,t} \leq I_{i,t} P_i^{\max}, \forall i,t \quad (8)$$

$$P_{i,t} - P_{i,(t-1)} \leq r_i^u(1 - y_{i,t}) + P_i^{\min} y_{i,t}, \forall i,t \quad (9)$$

$$-P_{i,t} + P_{i,(t-1)} \leq r_i^d(1 - z_{i,t}) + P_i^{\min} z_{i,t}, \forall i,t \quad (10)$$

$$(x_{i,t-1}^{\text{on}} - T_i^{\text{on}})(I_{i,t-1} - I_{i,t}) \geq 0, \forall i,t \quad (11)$$

$$(x_{i,t-1}^{\text{off}} - T_i^{\text{off}})(I_{i,t} - I_{i,t-1}) \geq 0, \forall i,t \quad (12)$$

$$0 \leq P_{r,t}^R \leq \tilde{P}_{r,t}^R, \forall r,t \Big\} \quad (13)$$

Equation (5) stands for the load balance constraint, where $d_{m,t}$ is the load demand at bus $m$ at time $t$. The net power injection is modeled in (6), where $\mathcal{G}_m$ and $\mathcal{R}_m$ denote the thermal unit set and RES unit set at bus $m$, respectively. The network constraints are enforced in (7). Equation (8) model the upper and lower power output limits of the thermal units. The thermal units are enforced by the sequential ramping rate limits in (9)-(10), where $I_{i,t}$, $y_{i,t}$, and $z_{i,t}$ are the indicators of the unit being on, started-up, and shut-down, respectively. Equations (9) and (10) show that a unit has to operate at its minimum capacity in two cases: right after it is turned on or right before it is turned off, which implies that the unit cannot provide reserve in those two cases. The minimum on/off time constraints are modeled in (11)-(12). The dispatchable RES output respects the constraint (13). The ramping rates of the RES output are not modeled in this paper. If the ramping constraints of RES similar to those of thermal units are enforced for RES, all the conclusions in this paper are still valid.

The feasible region $\mathcal{F}_2$ in (4) is defined as

$$\mathcal{F}_2 := \Big\{(P, I, P^R) : \forall \boldsymbol{\epsilon} \in \mathcal{U}, \exists \{\hat{P}_{i,t}, P^F_{i,t}, \hat{P}^R_{i,t}\} \text{ such that}$$

$$\sum_i \hat{P}_{i,t} + \sum_r \hat{P}^R_{r,t} + \sum_f P^F_{f,t} = \sum_m d_{m,t}, \forall t \quad (14)$$

$$\hat{P}^{\text{inj}}_{m,t} = \sum_{i \in \mathcal{G}_m} \hat{P}_{i,t} + \sum_{r \in \mathcal{R}_m} \hat{P}^R_{r,t} + \sum_{f \in \mathcal{Q}_m} P^F_{f,t} - d_{m,t}, \forall m \quad (15)$$

$$-F_l \leq \sum_m \Gamma_{l,m} \hat{P}^{\text{inj}}_{m,t} \leq F_l, \forall l, t. \quad (16)$$

$$I_{i,t} P^{\min}_i \leq \hat{P}_{i,t} \leq I_{i,t} P^{\max}_i, \forall i, t \quad (17)$$

$$-I_{i,t} R^d_i \leq \hat{P}_{i,t} - P_{i,t} \leq I_{i,t} R^u_i, \forall i, t \quad (18)$$

$$I^F_{f,t} P^{\min,F}_f \leq P^F_{f,t} \leq I_{f,t} P^{\max,F}_f, \forall f, t \quad (19)$$

$$0 \leq \hat{P}^R_{r,t} \leq \tilde{P}^R_{r,t} + \epsilon_{r,t}, \forall r, t \Big\}. \quad (20)$$

where $\hat{P}_{i,t}$ and $\hat{P}^R_{i,t}$ are the re-dispatches of thermal units and RES units, respectively, when the uncertainties are revealed. (14)-(15) are the load balance constraint, net power injection, and network constraints, respectively, after the uncertainty is revealed. The re-dispatch of thermal unit is limited by the ramping rate in (18). The $P^F_{f,t}$ denotes the power output of the fast startup units. It is enforced by the upper and lower generation limit in (19), where $I^F_{f,t}$ is the on/off indicator. The RES power level is not greater than the available maximum output in (20). It is noted that the ramping limit between $\hat{P}_{i,t}$ and $\hat{P}_{i,t+1}$ is not modeled. Instead, we model the ramping constraints from $P_{i,t}$ to $P_{i,t+1}$, from $P_{i,t}$ to $\hat{P}_{i,t}$, and from $P_{i,t+1}$ to $\hat{P}_{i,t+1}$. In the DAM, it is assumed that the operator has enough time to adjust the dispatches.

An advantage of the model formulated in this paper is that the deliverable ramping capability is guaranteed at each time interval. It becomes useful when considering the causality of the economic dispatch. In the existing two-stage robust approach, the causality is ignored. In comparison with the multi-stage robust approaches [19], the full recourse actions instead of AP policies are modeled in this paper which leads to more optimistic solution. It should be noted that the ramping constraint in this paper is less precise as ramping constraint between $\hat{P}_{i,t}$ and $\hat{P}_{i,t+1}$ is relaxed. This is consistent with the industry practice which is also called "corrective" dispatch after contingencies occur. It is worth mentioning that the solution approach in the next section is still applicable when the ramping constraint between $\hat{P}_{i,t}$ and $\hat{P}_{i,t+1}$ is enforced.

## III. SOLUTION APPROACH

In this section, we show that the two-stage RSCUC model (P1) can be converted into a single-level MIP problem. Accordingly, the solution to the original two-stage RSCUC problem can be obtained by solving the converted single-level MIP without using the Bender's decomposition or the Column Generation (CG) framework [4], [5], [13]. To solve (P1), a feasibility check problem regarding $\mathcal{F}_2$ is established as

$$\text{(FC)} \max_{\boldsymbol{\epsilon} \in \mathcal{U}} \min_{(\hat{P}, \hat{P}^R, P^F, s) \in \mathcal{F}_3(P, I, P^R, \boldsymbol{\epsilon})} \sum_t \sum_m s_{m,t}, \quad (21)$$

where $\mathcal{F}_3(P, I, P^R, \boldsymbol{\epsilon})$ is defined as

$$\mathcal{F}_3(P, I, P^R, \boldsymbol{\epsilon}) := \Big\{(\hat{P}, \hat{P}^R, P^F, s) :$$

$$\sum_i \hat{P}_{i,t} + \sum_r \hat{P}^R_{r,t} + \sum_f P^F_{f,t} = \sum_m d_{m,t} - s_{m,t}, \forall t \quad (22)$$

$$-F_l \leq \sum_m \Gamma_{l,m} (\hat{P}^{\text{inj}}_{m,t} + s_{m,t}) \leq F_l, \forall l, t \quad (23)$$

$$0 \leq \hat{P}^R_{r,t} \leq \tilde{P}^R_{r,t} + \epsilon_{r,t}, \forall r, t \quad (24)$$

$$s_{m,t} \geq 0, \forall m, t \quad (25)$$

$$(15), (17), (18), (19) \Big\}.$$

It is observed that the problem (FC) is a two-level problem. The out-level maximization problem is to find out the worst $\boldsymbol{\epsilon}$ which leads to the largest load curtailment (i.e. summation of $s_{m,t}$). The inner-level minimization problem is to find out the re-dispatch solution which leads to the lowest load curtailment. In comparison to the model shown in [17], no traditional generation curtailment (i.e. slack variable) is formulated. The reason is that the output of the dispatchable RES can be reduced and no over generation occurs. The relation between the problem (FC) and feasible region $\mathcal{F}_2$ can be described by the following proposition.

**Proposition 1.** $\{P_{i,t}, I_{i,t}, P^R_{r,t}\} \in \mathcal{F}_2$, if and only if the optimal value to the max-min problem (FC) is zero.

To solve the problem (FC), the inner problem is normally converted to its dual problem in literatures [4], [5], [17], where a bilinear problem is therefore formulated. In this paper, since the fast startup units are modeled in the re-dispatch process, the strong duality doesn't hold. Therefore, the traditional approach cannot be applied. Fortunately, the optimal solution to (FC) can be identified by the following theorem.

**Theorem 1.** The optimal solution to (FC) is obtained when $\epsilon_{r,t} = -u_{r,t}, \forall r, t$.

*Proof.* Denote the solution to (FC) as $\{\boldsymbol{\epsilon}^* : \epsilon^*_{r,t} = -u_{r,t}\}$, $z(P, I, P^R, \boldsymbol{\epsilon}^*)$ stands for the optimal value to the inner MIP problem when $\boldsymbol{\epsilon} = \boldsymbol{\epsilon}^*$. Assume there exists a point $\boldsymbol{\epsilon}' \neq \boldsymbol{\epsilon}^*$ where $\epsilon_{r',t} > -u_{r,t}$ for RES $r'$, the optimal value $z(P, I, P^R, \boldsymbol{\epsilon}') > z(P, I, P^R, \boldsymbol{\epsilon}^*)$. According to (24),

$$\{\hat{P}_{r',t} : 0 \leq \hat{P}^R_{r',t} \leq \tilde{P}^R_{r',t} - u_{r',t}\} \subset \{\hat{P}_{r',t} : 0 \leq \hat{P}^R_{r',t} \leq \tilde{P}^R_{r',t} + \epsilon_{r',t}\}.$$

As other constraints remain the same,

$$\mathcal{F}_3(P, I, P^R, \boldsymbol{\epsilon}^*) \subseteq \mathcal{F}_3(P, I, P^R, \boldsymbol{\epsilon}') \quad (26)$$

holds. The feasible region is enlarged, then the optimal values to the inner minimization problem have following relation

$$z(P, I, P^R, \boldsymbol{\epsilon}') \leq z(P, I, P^R, \boldsymbol{\epsilon}^*). \quad (27)$$

It can be observed that it is contradicted with the assumption. Hence, $\epsilon^*$ is the optimal solution to (FC), which leads to the largest load curtailment. □

The implication of Theorem 1 is that the worst case is fixed and is independent of the first stage solution $\{P_{i,t}, I_{i,t}, P^R_{r,t}\}$. Therefore, we can directly add the worst case into the RSCUC without solving the computation intensive max-min problem. Define

$$\mathcal{F}_4 := \Big\{(P, I, P^R, \hat{P}, \hat{P}^R, P^F) : 0 \leq \hat{P}^R_{r,t} \leq \tilde{P}^R_{r,t} - u_{r,t}$$
$$(14), (15), (16), (17), (18), (19)\Big\}.$$

Then the two-stage RSCUC (P1) is converted into a single-level MIP problem as follows.

$$\text{(P2):} \quad \min \quad \left\{ \begin{array}{l} \sum_i \sum_t \Big(C^P_i(P_{i,t}) + C^I_i(I_{i,t})\Big) \\ + \sum_r \sum_t C^R_r(P^R_{i,t}) \end{array} \right\}$$
$$\text{s.t.} \quad \{P_{i,t}, I_{i,t}, P^R_{i,t}\} \in \mathcal{F}_1,$$
$$\{P_{i,t}, I_{i,t}, P^R_{i,t}, \hat{P}_{i,t}, \hat{P}^R_{i,t}, P^F_{f,t}\} \in \mathcal{F}_4$$

If the cost of the worst case is considered similar to [14], then the new problem can be formulated as

$$\text{(P3):} \quad \min \quad (1-w) \left\{ \begin{array}{l} \sum_i \sum_t \Big(C^P_i(P_{i,t}) + C^I_i(I_{i,t})\Big) \\ + \sum_r \sum_t C^R_r(P^R_{i,t}) \end{array} \right\}$$
$$+ w \left\{ \begin{array}{l} \sum_i \sum_t \Big(C^P_i(\hat{P}_{i,t}) + C^I_i(I_{i,t})\Big) \\ + \sum_r \sum_t C^R_r(\hat{P}^R_{i,t}) + \sum_f \sum_t C^f_f(P^F_{f,t}) \end{array} \right\}$$
$$\text{s.t.} \quad \{P_{i,t}, I_{i,t}, P^R_{i,t}\} \in \mathcal{F}_1,$$
$$\{P_{i,t}, I_{i,t}, P^R_{i,t}, \hat{P}_{i,t}, \hat{P}^R_{i,t}, P^F_{f,t}\} \in \mathcal{F}_4.$$

In comparison to the CG based framework, the main difference is that only one additional scenario is considered, which guarantees the robustness (i.e. system can survive when the uncertainty is revealed without load curtailment). Therefore, the increased computation burden due to the robustness is acceptable.

As the RESs are dispatchable, they can also be used to accommodate the uncertainties. Hence, RESs do not require FRs any more, they actually function as FRs. This is especially critical for systems with high level RESs.

## IV. CASE STUDY

The proposed RSCUC model and solution approaches are tested with the modified IEEE 118-bus system. The system is consisted of 54 thermal units, 186 branches, 15 wind farms, and 15 solar farms. The MIP solver Gurobi 5.6.3 [25] is utilized to solve the MIP problems on PC with Intel i7-3770@3.40GHz 8GB RAM.

The peak load is 6600MW in 24 hours for the modified IEEE 118-bus system. Hour 11 and Hour 20 are two peak points. Both the solar farms and the wind farms are simulated simultaneously in this section. The aggregated wind power output forecast is obtained from ISO-NE. The data in scenario 1 and scenario 2 is one fifth of the day-ahead forecasts for wind power output on 06-17-2015 and 03-31-2015, respectively. The solar farm output is from NREL. The aggregated solar power is the output from Arizona on 01-01-2006. The 15 wind farms and 15 solar farms are located at 30 different buses. The detailed data including generating unit parameters, line reactance and ratings, and load profiles can be found at http://motor.ece.iit.edu/Data/ROSCUC_118.xls.

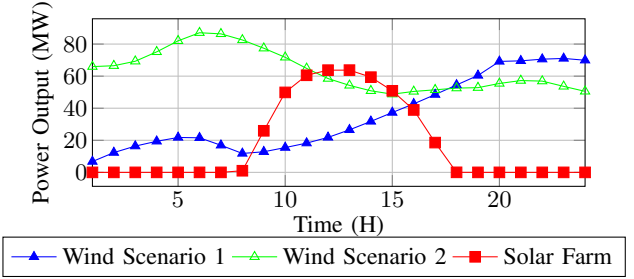

Fig. 4. Aggregated Wind and Solar Power Output Forecasting

The curves in Figure 4 depict the expected value of the output $\tilde{P}^R_{r,t}$. The actual available RES power is $\tilde{P}^R_{r,t} + \epsilon_{r,t}$. In this section, the uncertainty $\epsilon_{r,t}$ is assumed within the confidence interval, i.e. $-\alpha \tilde{P}^R_{r,t} \leq \epsilon_{r,t} \leq \alpha \tilde{P}^R_{r,t}$. $\alpha \in \mathbb{R}$ is chosen from 20% to 30% in the simulations. Denote $\beta$ as the RES level. $\beta = 1$ stands for the base RES level.

### A. Dispatchable RES V.S. Non-Dispatchable RES

In many robust literatures, the generation outputs from RES are modeled as inelastic load [4], [5], [15], [16]. Any fluctuation from RES output must be accommodated with flexible resources. In contrast, the RESs are assumed to be dispatchable within their available power level in this paper. In this part, the main objective is to compare the impacts of the RES dispatches. For simplicity, the fast startup generators are not considered in the second stage, which may cause computation issue for the traditional RSCUC. At the same time, the bids of the RESs are assume to be zeros. The confidence interval is set to 25%, and $w = 0$.

The simulations are performed in two steps. In the first step, the Robust UC and Dispatch are determined based on (P2) or other robust approaches with non-dispatchable RESs [5], [4], [17]. In the second step, 1000 samples of $\epsilon$ are generated for the RESs following the normal distribution. Then the available flexibilities are re-dispatched with the realized RES outputs. We perform the sensitivity analysis with the RESs output level in terms of operation cost and UC hours.

Table I presents the base cost and average cost comparisons between solutions with dispatchable and non-dispatchable RESs with wind scenario 1. The "base cost" is the total operation cost to supply the load with the expected RES output levels. The "average cost" is the average total operation costs for all the sample points. Column "$\beta$" shows the increasing RES level. The total expected available power from wind and solar is 19,197 MWh, which is 14.6% of the total load demand 131,472MWh. When $\beta$ reaches 5, the available power from wind and solar is 73% of the total load demand. It can be observed that the cost for the RSCUC with dispatchable RES in this paper is always smaller than the one with non-dispatchable RES. It is also consistent with the feasible region comparison shown in Figure 3. With smaller feasible unit dispatch region, the traditional RSCUC normally is more conservative. The total UC hours ("UCs(h)") represents the total committed hours for all units. In order to provide more flexible resources, the traditional RSCUC with non-dispatchable RESs has to keep more units ON. In comparison, the total UC hours of the approach in this paper are always smaller.

TABLE I
THE BASE CASE COST AND AVERAGE COST COMPARISON WITH
INCREASING RES OUTPUTS (SCENARIO 1)

| $\beta$ | Base Cost ($) | | Average Cost ($) | | UCs (h) | |
|---|---|---|---|---|---|---|
| | Dis. | Non-Dis. | Dis. | Non-Dis. | Dis. | Non-Dis. |
| 1.0 | 1,613,353 | 1,613,397 | 1,570,180 | 1,570,617 | 625 | 641 |
| 1.6 | 1,415,987 | 1,419,017 | 1,380,330 | 1,381,957 | 629 | 642 |
| 1.7 | 1,386,416 | 1,400,034 | 1,351,274 | 1,363,356 | 644 | 704 |
| 2.2 | 1,253,351 | NaN | 1,212,849 | NaN | 658 | NaN |
| 5.0 | 697,865 | NaN | 663,279 | NaN | 388 | NaN |

An important observation is that the cost differences change with the RES production levels. For the base level (i.e. $\beta = 1$), the base-case cost difference is $1,613,397 - $1,613,353 = $44, and the average cost difference is $1,570,617 - $1,570,180 = $437. When the RES production level increases to 1.7, the solution with dispatchable RES can save $1,400,034 - $1,386,416=$13,618 in base-case cost in comparison to the one with non-dispatchable RES. The average cost saving increases to $1,351,274 - $1,363,356=$12,082 from $437. It demonstrate that the proposed approach with dispatchable RES leads much lower cost when the RES level is high.

Another observation is that the RSCUC with non-dispatchable RES becomes infeasible when the RES production level increases to 2.2. It is assumed that the confidence interval $\epsilon_{r,t}$ is [-25%, 0.25%] of the available production level. Hence, with larger expected power output, the uncertainties the system should accommodate also increase. However, the ramping capability, which is provided by the thermal unit, has an upper limit. When the uncertainty is above the limit, no feasible solution can be found for the RSCUC with non-dispatchable RES. The system operators may have to curtail the load or shut down thermal units regardless of the UC schedules. In contrast, the problem with dispatchable RES is feasible even when the RES level reaches 5.0. The reason is that the RESs can also provide flexibilities. At a level of 5.0, the total average cost is reduced to $663,279, which is 42.24% of the base-case average cost $1,570,180.

Table II shows another set of cost results with wind scenario 2 depicted Figure 4. The total expected available RES production is 22.3% of the load demand. It can be observed that when RES level is over 1.2, the RSCUC with non-dispatchable RES becomes infeasible. When the RES level increases to 5.0, the total expected available RES power is 111.5% of the load demand. There are only 248 UC hours in this case and the total cost drops to $324,916. Although the RES power is larger than the load demand, two factors prevent the cost from being zero. One is that there are no energy storages to store the over-generated power. For example, wind power reaches its peak point while the load demand is lower at 6:00AM. The other one is the line congestions, which also prevent free RES energy from all being delivered to certain load buses.

TABLE II
THE BASE-CASE COST AND AVERAGE COST COMPARISON WITH
INCREASING RES OUTPUTS (SCENARIO 2)

| $\beta$ | Base Cost ($) | | Average Cost ($) | | UCs (h) | |
|---|---|---|---|---|---|---|
| | Disp. | Non-Disp. | Disp. | Non-Disp. | Disp. | Non-Disp. |
| 1.0 | 1,446,969 | 1,447,583 | 1,410,778 | 1,411,577 | 614 | 631 |
| 1.2 | 1,352,607 | 1,357,410 | 1,317,716 | 1,322,519 | 606 | 638 |
| 1.3 | 1,308,989 | NaN | 1,274,367 | NaN | 592 | NaN |
| 5.0 | 351,987 | NaN | 324,916 | NaN | 248 | NaN |

TABLE III
PROCURED RES ENERGY WITH INCREASING RES LEVELS AND
DECREASING RES BIDS ($\alpha = 25\%, w = 0$)

| Bid[a] | $\beta = 1.0$ | | | $\beta = 2.0$ | | | $\beta = 5.0$ | | |
|---|---|---|---|---|---|---|---|---|---|
| | RES[b] | (%) | UCs[c] | RES | (%) | UCs (h) | RES | (%) | UCs |
| 5 | 29,321 | 100.00 | 614 | 53,908 | 91.93 | 554 | 103,824 | 70.82 | 213 |
| 0 | 29,321 | 100.00 | 614 | 54,542 | 93.01 | 635 | 104,364 | 71.19 | 248 |
| -5 | 29,321 | 100.00 | 614 | 54,583 | 93.08 | 640 | 105,270 | 71.81 | 304 |

[a] $/MWh; [b] MWh; [c] h

### B. Impacts of the RES Bids

In the DAM, the energies are traded at the financially binding prices between market participants. In this part, the bids of the RES are simulated to show differences of the energy settlements. The simulations are analyzed and discussed from RES energy accommodation point of view rather than the market cost perspective. Hence, three different simple bid offers (i.e. $5/MW, $0/MW, and $-5/MW) are simulated, and the typical wind scenario 2 is used. As an example of negative bidding price, the ISO New England allows the market participants to bid at negative prices since December 3, 2014. Four groups of sensitivity analysis for RES procurement are discussed in this part:

- with increasing RES levels;
- with increasing uncertainty levels;
- with increasing weight factor of the worst case;
- with increasing fast startup units.

*1) Increasing RES Levels:* The procured RES energy and UC hours With increasing RES level are presented in Table III. The RES levels are chosen as 1.0, 2.0, and 5.0. The procured RES energy in DAM are shown in the columns denoted as "RES". They are the same with different bid offers when the production level $\beta$=1. However, if the RES increases to 2.0 times of the base level, then different bid offers result in different procured RES energy. Column "%" shows the procured RES energy as a percentage of the expected available RES energy. For example, the procured RES energy is 53,908 MWh with a $5/MWh bid offer, which is 91.93% of the

TABLE IV
PROCURED RES ENERGY WITH INCREASING UNCERTAINTY
LEVELS AND DECREASING RES BIDS ($\beta = 2.0, w = 0$)

| Bid[a] | $\alpha$=0% | | | $\alpha$=20% | | | $\alpha$=30% | | |
|---|---|---|---|---|---|---|---|---|---|
| | RES[b] | (%) | UCs[c] | RES | (%) | UCs | RES | (%) | UCs |
| 5 | 57,867 | 98.68 | 388 | 55,524 | 94.68 | 457 | 52,071 | 88.79 | 633 |
| 0 | 57,909 | 98.75 | 388 | 56,423 | 96.22 | 565 | 52,293 | 89.17 | 661 |
| -5 | 57,971 | 98.86 | 388 | 56,626 | 96.56 | 589 | 52,390 | 89.34 | 674 |

[a] $/MWh;  [b] MWh;  [c] h.

TABLE V
PROCURED RES ENERGY WITH INCREASING WEIGHT FACTORS FOR THE
WORST-CASE COST AND DECREASING RES BIDS ($\alpha = 25\%, \beta = 2.0$)

| Bid[a] | $w = 0$ | | | $w = 0.2$ | | | $w = 1$ | | |
|---|---|---|---|---|---|---|---|---|---|
| | RES[b] | (%) | UCs[c] | RES | (%) | UCs | RES | (%) | UCs |
| 5 | 53,908 | 91.93 | 554 | 53,399 | 91.06 | 481 | 47,010 | 80.16 | 421 |
| 0 | 54,542 | 93.01 | 635 | 54,040 | 92.15 | 562 | 47,010 | 80.16 | 421 |
| -5 | 54,583 | 93.08 | 640 | 54,542 | 93.01 | 623 | 47,010 | 80.16 | 421 |

[a] $/MWh;  [b] MWh;  [c] h.

TABLE VI
PROCURED RES ENERGY WITH INCREASING NUMBER OF FAST STARTUP
UNITS AND DECREASING RES BIDS ($\alpha = 25\%, \beta = 2.0, w = 0.2$)

| Bid[a] | # of FU=0 | | | # FU=3 | | | # of FU=6 | | |
|---|---|---|---|---|---|---|---|---|---|
| | RES[b] | (%) | UCs[c] | RES | (%) | UCs | RES | (%) | UCs |
| 5 | 53,399 | 91.06 | 481 | 56,918 | 97.06 | 411 | 57,799 | 98.56 | 394 |
| 0 | 54,040 | 92.15 | 562 | 57,018 | 97.23 | 416 | 57,804 | 98.57 | 394 |
| -5 | 54,542 | 93.01 | 623 | 57,427 | 97.93 | 465 | 57,842 | 98.64 | 393 |

[a] $/MWh;  [b] MWh;  [c] h.

expected available RES energy. In contrast, if the bid offer decreases to $0/MWh, then the procured RES energy increases by 54,542 MWh - 53,908 MWh = 634 MWh. In other words, the thermal units supply 634 MWh fewer load demands. Table III also shows that the UC hours increase to 635 from 554 in this case. An interesting observation is that fewer thermal units are committed to supply more loads when the bid offer changes from $0/MWh to $5/MWh. The reason is that, to accommodate more cheaper $0/MWh RES energy, more generation reserves must be kept in order to accommodate the uncertainties. Thus, more units are committed, which can provide more ramping capabilities. It indicates that the opportunity cost of keeping these reserves is smaller than that of curtailing the RES energies. A similar trend can also be observed when $\beta$=5.0.

With increasing $\beta$, the procured percentage of expected available RES energy is decreasing. For example, when $\beta = 1$, all RES energy is procured in DAM. In contrast, when $\beta = 5.0$, only around 71% of the available RES energy is scheduled to be delivered into the system. It indicates that when the RES level is high, the most economical way is to have some RES energy spilled. It becomes expensive to accommodate all the RESs as the reserve cost is high at this time. It is worth mentioning that the total operation cost is negative when $\beta$=5.0 and the bid offer is $-5/MWh.

*2) Increasing Uncertainty Levels:* As mentioned above, more generation reserves are required if more energy from uncertain RESs are delivered into the grid. In this paragraph, we discuss the impacts of the uncertainty levels on RES energy procurement in DAM. Table IV presents the procured RES energy in DAM with different bid offers and uncertainty levels. It is assumed that the RES level $\beta = 2.0$, which is about 44.6% of the total load demand. It can be observed with the same bid offer, the increase of the uncertainty level causes the drop of the amount of procured RES energy. For example, with $0/MWh bid offer, the procured RES energy is 56,423 MWh, and total UC hours are 565 when $\alpha = 20\%$. When $\alpha$ increases to 30%, then the procured RES energy drops to 52,293 MWh, and the total UC hours increase to 661. Similar trends are also observed when the bid offer is $5/MWh or $-5/MWh. In general, more uncertainties are accommodated by more ramping capabilities in the system.

Another observation is that with the increasing uncertainty, the procured RES energy percentage (%) is decreasing. A higher uncertainty level requires more spinning reserves in the system to ensure the solution robustness. With the assumption $\alpha = 0\%$, there are no uncertainties for the RES energy. The procured RES energy percentage is 98.75% with zero bid offer.

When the uncertainty level $\alpha = 30\%$, the procured percentage drops by 9.58% to 89.17%. It indicates that compared with the cost of the supporting reserve for RES energy, it is much cheaper to schedule the thermal units to supply the load. In Figure 5, the spillage of the RES energy is depicted. It shows that the RES energy spillages occur from 1:00AM to 17:00AM. At Hour 1, the ON/OFF statuses of many units are changed from those at Hour 0. With the assumption that the spinning reserve is zero when the unit is at its first ON hour or the last OFF hour, the reserves provided at Hour 1 is small. Consequently, the spillage is the largest at Hour 1.

*3) Increasing Weight Factors for the Worst-Case Cost:* In Table III and Table IV, the weight factor for the worst-case cost is set to 0. In this part, we perform the sensitivity analysis of RES energy procurement with the increasing weight factors for the worst-case cost. The simulation results are presented in Table V. In general, the increase of the weight factor results in fewer RES energy procurements. When $w = 0$, only the base case is considered. For example, with bid offer at $5/MWh, the procured RES energy percentage decreases from 91.93% to 80.16% when $w$ is changed from 0 to 1. In the extreme case $w = 1$, only the worst case is optimized. No upward supporting reserves for RES energy are needed in this case. Therefore, we have fewer UC hours. Consequently, the system can accommodate fewer RES energies in the base case. Similar trends can also be observed when the bid offers are $0/MWh and $-5/MWh.

*4) Increasing Fast Startup Units:* When considering the fast startup units, the system has more flexible resources, although the price of these resources is expensive. In Table VI, the RES energy procured in the DAM is presented with the increasing number of fast startup units. It can be observed that the increase of the fast startup units leads to more RES energy procurements. For example, the procured RES energy percentage increases to 98.57% from 91.06% with $5/MWh bid offer when the number of fast startup units changes to 6 from 0. On the other hand, the UC hours are decreasing. It indicates that more flexible resources help the system

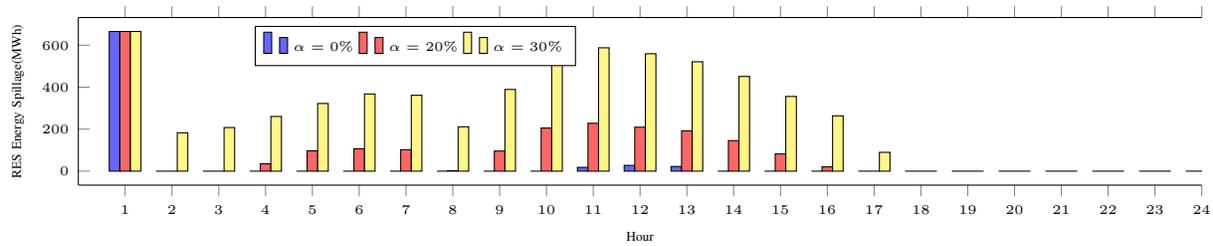

Fig. 5. RES Energy Spillage w.r.t. Different Uncertainty Level ($\beta = 2.0, w = 0$, bid=$/MWh).

accommodate more RES energies.

*C. Computation Comparison*

The robustness of the approach in this paper can be guaranteed by an extra scenario (i.e. the worst case). With the default parameters, the MIP solver can get the solution within 2 seconds for the IEEE 118-Bus system. The existing RSCUC with non-dispatchable RES normally costs more time to get the solution. Without any acceleration techniques, approaches in [4], [13], [5] cannot get the solution in 2 hours. Two factors slow down the solution process in the these approaches. The first one is that these algorithms often involve iterative processes. For example, in the Benders Decomposition approach [5], [4], [23], a large number of new MIP problems with Benders cuts generated by solving subproblems must be solved repeatedly. In the Column Generation approach [13], [17], different worst points are calculated each time. Several increasing-size MIP problems need to be solved repeatedly in the iterative process. The second one is that solving the non-convex max-min subproblem is a challenging task. In the approach proposed in this paper, no max-min subproblem needs to be solved.

## V. CONCLUSIONS

In this paper, a robust SCUC and dispatch model considering the RES bids are proposed in the system with high RES level. It is proved that with the dispatchable RES, the worst case for the second stage can be directly identified. The conclusion always holds whether the strong duality holds for the re-dispatch problem or not. With this conclusion, the robustness can be guaranteed by adding only one extra scenario in the original SCUC problem.

The simulation results in this paper show that by increasing the RES level, the total cost can be lowered. However, there is an increasing chance that the traditional robust SCUC with non-dispatchable RES is infeasible if the RES level is high. The simulation results also show that it is not the most economic strategy to accommodate all the power from RESs into the grid. By taking the advantage of the dispatchable RES, this paper shows that the RES can also provide flexibilities. In the future power system with high RES penetration, the dispatchability of the RES will play a crucial role. It is a promising and attractive new research topic.